\documentclass{amsart}
\usepackage{amssymb,dbnsymb,multicol,graphicx,picins}
\usepackage[all]{xy}

\newlength{\standardunitlength}
\catcode`\@=11
\long\def\@makecaption#1#2{%
    \vskip 10pt
    \setbox\@tempboxa\hbox{%\ifvoid\tinybox\else\box\tinybox\fi
      \small\sf{\bfcaptionfont #1. }\ignorespaces #2}%
    \ifdim \wd\@tempboxa >\captionwidth {%
        \rightskip=\@captionmargin\leftskip=\@captionmargin
        \unhbox\@tempboxa\par}%
      \else
        \hbox to\hsize{\hfil\box\@tempboxa\hfil}%
    \fi}
\font\bfcaptionfont=cmssbx10 scaled \magstephalf
\newdimen\@captionmargin\@captionmargin=2\parindent
\newdimen\captionwidth\captionwidth=\hsize
\catcode`\@=12

\newlength{\globalparindent}
\newcommand{\udot}{{\mathaccent\cdot\cup}}

\theoremstyle{plain}
\newtheorem{theorem}{Theorem}
\newtheorem{conjecture}[theorem]{Conjecture}
\newtheorem{problem}[theorem]{Problem}

\theoremstyle{definition}
\newtheorem{definition}[theorem]{Definition}

\def\bbQ{{\mathbb Q}}
\def\bbR{{\mathbb R}}
\def\calA{{\mathcal A}}
\def\calB{{\mathcal B}}
\def\calD{{\mathcal D}}
\def\calE{{\mathcal E}}
\def\calK{{\mathcal K}}
\def\calP{{\mathcal P}}
\def\calS{{\mathcal S}}
\def\calT{{\mathcal T}}
\def\calU{{\mathcal U}}
\def\frakg{{\mathfrak g}}

\newcommand{\eps}[2]{{\hspace{-3pt}\begin{array}{c}%
  \raisebox{-2.5pt}{\includegraphics[#1]{figs/#2.eps}}%
\end{array}\hspace{-3pt}}}

\def\arXiv#1{arXiv:\linebreak[0]#1}

\def\sumint{\setbox0=\hbox{$\displaystyle\sum$}\mathop{\rlap{\copy0}
\kern0pt \hbox to \wd 0{\hss$\displaystyle\int$\hss}}}

\begin{document}

\title{Finite Type Invariants}

\author{Dror~Bar-Natan}
\address{
  Department of Mathematics\\
  University of Toronto\\
  Toronto Ontario M5S 3G3\\
  Canada
}
\email{drorbn@math.toronto.edu}
\urladdr{http://www.math.toronto.edu/\~{}drorbn}

\date{
  First edition: Aug.~13,~2004.
  This edition: \today
}

\thanks{This work was partially supported by NSERC grant RGPIN 262178.
  Electronic versions: {\tt
  http://\linebreak[0]www.math.toronto.edu/\linebreak[0]$\sim$drorbn/\linebreak[0]papers/\linebreak[0]EMP/} and arXiv:math.GT/0408182.
  Copyright is retained by the author.
}

\begin{abstract}
This is an overview article on finite type invariants, written for the {\em
Encyclopedia of Mathematical Physics}.
\end{abstract}

\maketitle

\begin{multicols}{2}

\section{Introduction} \label{sec:intro}

Knots belong to sailors and climbers and upon further reflection,
perhaps also to geometers, topologists or combinatorialists.
Surprisingly, throughout the 1980s it became apparent that knots are
also closely related to several other branches of mathematics in
general and mathematical physics in particular. Many of these
connections (though not all!) factor through the notion of ``finite
type invariants'' (aka ``Vassiliev'' or ``Goussarov-Vassiliev''
invariants)~\cite{Goussarov:New, Goussarov:nEquivalence,
Vassiliev:CohKnot, Vassiliev:Book, BirmanLin:Vassiliev,
Kontsevich:Vassiliev, Bar-Natan:OnVassiliev}.

Let $V$ be an arbitrary invariant of oriented knots in oriented space
with values in some Abelian group $A$. Extend $V$ to be an invariant of
$1$-singular knots, knots that may have a single singularity that locally
looks like a double point $\doublepoint$, using the formula
\begin{equation} \label{eq:doublepoint}
  V(\doublepoint)=V(\overcrossing)-V(\undercrossing).
\end{equation}
Further extend $V$ to the set $\calK^m$ of $m$-singular knots (knots
with $m$ double points) by repeatedly using~\eqref{eq:doublepoint}.

\begin{definition} We say that $V$ is of type $m$ if its extension
$\left.V\right|_{\calK^{m+1}}$ to $(m+1)$-singular knots vanishes
identically. We say that $V$ is of finite type if it is of type $m$ for
some $m$.
\end{definition}

Repeated differences are similar to repeated derivatives and hence it
is fair to think of the definition of $\left.V\right|_{\calK^m}$ as
repeated differentiation. With this in mind, the above definition
imitates the definition of polynomials of degree $m$. Hence finite type
invariants can be thought of as ``polynomials'' on the space of knots.

As we shall see below, finite type invariants are plenty and powerful
and they carry a rich algebraic structure and are deeply related to Lie
algebras (Section~\ref{sec:Basic}). There are several constructions for
a ``universal finite type invariant'' and those are related to
conformal field theory, the Chern-Simons-Witten topological quantum
field theory and Drinfel'd's theory of associators and quasi-Hopf
algebras (Section~\ref{sec:Fundamental}). Finite type invariants have
been studied extensively (Section~\ref{sec:FurtherDirections}) and
generalized in several directions (Section~\ref{sec:BeyondKnots}). But
the first question on finite type invariants remains unanswered:

\begin{problem} Honest polynomials are dense in the space of functions. Are
finite type invariants dense within the space of all knot invariants? Do
they separate knots?
\end{problem}

In a similar way one may define finite type invariants of framed knots
(and ask the same questions).

\subsection{Acknowledgement} I wish to thank O.~Dasbach for a correction.

\section{Basic facts} \label{sec:Basic}

\subsection{Classical knot polynomials} The first (non trivial!) thing to
notice is that there are plenty of finite type invariants and they are at
least as powerful as all the standard knot polynomials\footnote{Finite type
invariants are like polynomials on the space of knots; the standard phrase
``knot polynomials'' refers to a different thing --- knot invariants with
polynomial values.} combined:

\begin{theorem} \cite{Bar-Natan:OnVassiliev, BirmanLin:Vassiliev}
Let $J(K)(q)$ be the Jones polynomial of a knot $K$ (it is
a Laurent polynomial in a variable $q$). Consider the power series
expansion $J(K)(e^x)=\sum_{m=0}^\infty V_m(K)x^m$. Then each coefficient
$V_m(K)$ is a finite type knot invariant. (And thus the Jones polynomial
can be reconstructed from finite type information).
\end{theorem}

A similar theorem holds for the Alexander-Conway, HOMFLY-PT and
Kauffman polynomials~\cite{Bar-Natan:OnVassiliev}, and indeed, for
arbitrary Reshetikhin-Turaev invariants~\cite{ReshetikhinTuraev:Ribbon,
Lin:QuantumGroups}. Though it is still unknown if the signature of a knot can
be expressed in terms of its finite type invariants.

\subsection{Chord diagrams and the Fundamental Theorem}
\label{subsec:ChordDiagrams} The top derivatives of a multi-variable
polynomial form a system of constants that determine that polynomial up
to polynomials of lower degree.  Likewise the $m$th derivative
$V^{(m)}:=V(\doublepoint\overset{m}{\cdots}\doublepoint)$ of a type $m$
invariant $V$ is a constant (for
$V(\doublepoint\overset{m}{\cdots}\doublepoint\overcrossing)
-V(\doublepoint\overset{m}{\cdots}\doublepoint\undercrossing)
=V(\doublepoint\overset{m+1}{\cdots}\doublepoint)=0$ so $V^{(m)}$ is
blind to 3D topology) and likewise $V^{(m)}$ determines $V$ up to
invariants of lower type. Hence a primary tool in the study of finite
type invariants is the study of the ``top derivative'' $V^{(m)}$, also
known as ``the weight system of $V$''.

\parpic[r]{$\eps{width=4cm}{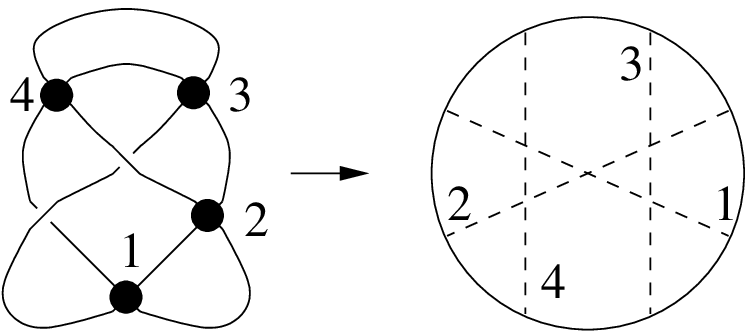}$}
Blind to 3D topology, $V^{(m)}$ only sees the combinatorics of the
circle that parametrizes an $m$-singular knot. On this circle there are $m$
pairs of points that are pairwise identified in the image; standardly one
indicates those by drawing a circle with $m$ chords marked (an ``$m$-chord
diagram'') as on the right.

\begin{definition} Let $\calD_m$ denote the space of all formal linear
combinations with rational coefficients of $m$-chord diagrams. Let
$\calA^r_m$ be the quotient of $\calD_m$ by all $4T$ and $FI$ relations
as drawn below (full details in e.g.~\cite{Bar-Natan:OnVassiliev}), and let
$\hat\calA^r$ be the graded completion of $\calA:=\bigoplus_m\calA^r_m$. Let
$\calA_m$, $\calA$ and $\hat\calA$ be the same as $\calA^r_m$,
$\calA^r$ and $\hat\calA^r$ but without imposing the $FI$ relations.
\[ 4T:\ \eps{height=8mm}{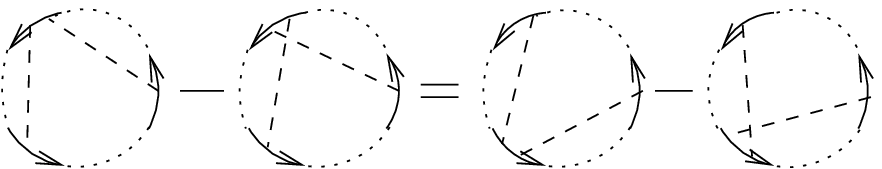}\quad FI:\ \eps{height=8mm}{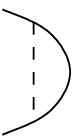} \]
\end{definition}

\begin{theorem} \label{thm:Fundamental} (The Fundamental Theorem)

\noindent $\bullet$\ (Easy part, \cite{Vassiliev:CohKnot, Goussarov:New,
BirmanLin:Vassiliev}). If $V$ is a rational valued type $m$ invariant
then $V^{(m)}$ defines a linear functional on $\calA^r_m$. If in
addition $V^{(m)}\equiv 0$, then $V$ is of type $m-1$.

\noindent $\bullet$\  (Hard part, \cite{Kontsevich:Vassiliev} and
Section~\ref{sec:Fundamental}). For any linear functional $W$ on
$\calA^r_m$ there is a rational valued type $m$ invariant $V$ so that
$V^{(m)}=W$.
\end{theorem}

Thus to a large extent the study of finite type invariants is reduced
to the finite (though super exponential in $m$) algebraic study of
$\calA^r_m$. A similar theorem reduces the study of finite type
invariants of framed knots to the study of $\calA_m$.

\subsection{The structure of $\calA$} Knots can be multiplied (the
``connected sum'' operation) and knot invariants can be multiplied.
This structure interacts well with finite type invariants and induces
the following structure on $\calA^r$ and $\calA$:

\begin{theorem} \cite{Kontsevich:Vassiliev, Bar-Natan:OnVassiliev,
Willerton:Hopf, ChmutovDuzhinLando:VasI} $\calA^r$ and $\calA$
are commutative and cocommutative graded bialgebras (i.e., each carries
a commutative product and a compatible cocommutative coproduct).
Thus both $\calA^r$ and $\calA$ are graded polynomial
algebras over their spaces of primitives, $\calP^r=\oplus_m\calP^r_m$ and
$\calP=\oplus_m\calP_m$.
\end{theorem}

Framed knots differ from knots only by a single integer parameter (the
``self linking'', itself a type $1$ invariant). Thus $\calP^r$ and $\calP$
are also closely related.

\begin{theorem} \cite{Bar-Natan:OnVassiliev}
$\calP=\calP^r\oplus\langle\theta\rangle$, where $\theta$ is the unique
$1$-chord diagram $\eps{width=5mm}{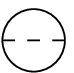}$.
\end{theorem}

\subsection{Bounds and computational results} The following table (taken
from~\cite{Bar-Natan:OnVassiliev, Kneissler:Twelve}) shows the number of
type $m$ invariants of knots and framed knots modulo type $m-1$ invariants
($\dim\calA_m^r$ and $\dim\calA_m$) and the number of multiplicative
generators of the algebra $\calA$ in degree $m$ ($\dim\calP_m$) for
$m\leq 12$. Some further tabulated results are
in~\cite{Bar-Natan:Computations}.

\begin{center}
{\small
  \def\n#1{{$\!\!#1\!\!$}}
  \begin{tabular}{||c|c|c|c|c|c|c|c|c|c|c|c|c|c||}
    \hline
    \n{m} &
	\n{0} & \n{1} & \n{2} & \n{3} & \n{4} & \n{5} & \n{6} & \n{7} &
	\n{8} & \n{9} & \n{10} & \n{11} & \n{12} \\
    \hline
    \n{\dim\calA_m^r} &
	\n{1} & \n{0} & \n{1} & \n{1} & \n{3} & \n{4} & \n{9} & \n{14} &
	\n{27} & \n{44} & \n{80} & \n{132} & \n{232} \\
    \n{\dim\calA_m} &
	\n{1} & \n{1} & \n{2} & \n{3} & \n{6} & \n{10} & \n{19} & \n{33} &
	\n{60} & \n{104} & \n{184} & \n{316} & \n{548} \\
    \n{\dim\calP_m} &
	\n{0} & \n{1} & \n{1} & \n{1} & \n{2} & \n{3} & \n{5} & \n{8} &
	\n{12} & \n{18} & \n{27} & \n{39} & \n{55} \\
    \hline
  \end{tabular}
}
\end{center}

Little is known about these dimensions for large $m$. There is an
explicit conjecture in~\cite{Broadhurst:ConjecturedEnumeration} but no
progress has been made in the direction of proving or disproving it.
The best asymptotic bounds available are:

\begin{theorem} For large $m$, $\dim\calP_m>e^{c\sqrt m}$ (for any fixed
$c<\pi\sqrt{\frac23}$)
\cite{Dasbach:CombinatorialStructureIII, Kontsevich:Unpublished} and
$\dim\calA_m<6^mm!\sqrt{m}/\pi^{2m}$
\cite{Stoimenow:Enumeration, Zagier:StrangeIdentity}.
\end{theorem}

\subsection{Jacobi diagrams and the relation with Lie algebras}
\label{subsec:Jacobi} Much of the richness of finite type invariants
stems from their relationship with Lie algebras.
Theorem~\ref{thm:Jacobi} below suggests this relationship on an
abstract level, Theorem~\ref{thm:Lie} makes that relationship
concrete and Theorem~\ref{thm:PBW} make is a bit deeper.

\begin{theorem} \label{thm:Jacobi} \cite{Bar-Natan:OnVassiliev}
The algebra $\calA$ is isomorphic to the
algebra $\calA^t$ generated by ``Jacobi diagrams in a circle'' (chord
diagrams that are also allowed to have oriented internal trivalent
vertices) modulo the $AS$, $STU$ and $IHX$ relations. See the figure
below.
\end{theorem}

\vskip -1mm
\[ \eps{width=2.8in}{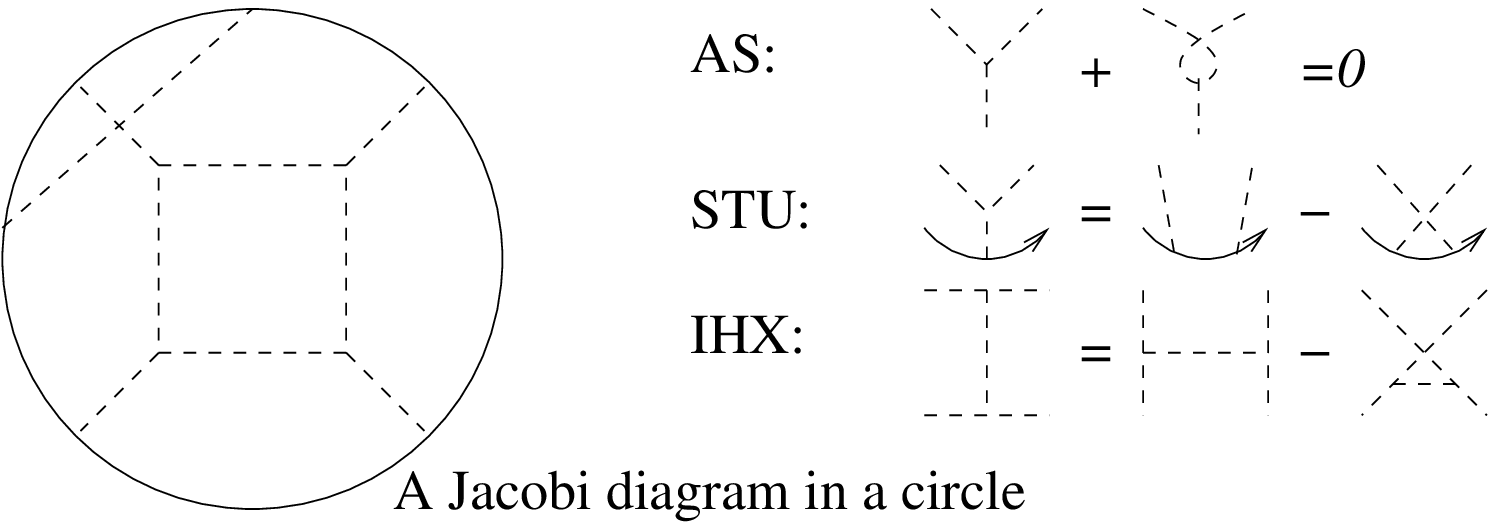} \]
\vskip 1mm

Thinking of trivalent vertices as graphical analogs of the Lie bracket, the
$AS$ relation become the anti-commutativity of the bracket, $STU$
become the equation $[x,y]=xy-yx$ and $IHX$ becomes the Jacobi identity.
This analogy is made concrete within the proof of the following:

\begin{theorem} \label{thm:Lie} \cite{Bar-Natan:OnVassiliev}
Given a finite dimensional metrized Lie
algebra $\frakg$ (e.g., any semi-simple Lie algebra) there is a map
$\calT_\frakg:\calA\to\calU(\frakg)^\frakg$ defined on $\calA$ and
taking values in the invariant part $\calU(\frakg)^\frakg$ of the
universal enveloping algebra $\calU(\frakg)$ of $\frakg$. Given also a
finite dimensional representation $R$ of $\frakg$ there is a linear
functional $W_{\frakg,R}:\calA\to\bbQ$.
\end{theorem}

The last assertion along with Theorem~\ref{thm:Fundamental} show that
associated with any $\frakg$, $R$ and $m$ there is a weight system and
hence a knot invariant. Thus knots are unexpectedly linked with Lie
algebras.

The hope \cite{Bar-Natan:OnVassiliev} that all finite type invariants
arise in this way was dashed by~\cite{Vogel:Structures,
Vogel:UniversalLieAlgebra, Lieberum:NotComing}. But finite type invariants
that do not arise in this way remain rare and not well understood.

The Poincar\'e-Birkhoff-Witt (PBW) theorem of the theory of Lie
algebras says that the obvious ``symmetrization'' map
$\chi_\frakg:\calS(\frakg)\to\calU(\frakg)$ from the symmetric algebra
$\calS(\frakg)$ of a Lie algebra $\frakg$ to its universal enveloping
algebra $\calU(\frakg)$ is a $\frakg$-module isomorphism. The following
definition and theorem form a diagrammatic counterpart of this theorem:

\parpic(12mm,13mm)[r]{\raisebox{-14mm}{$
  \eps{width=11mm}{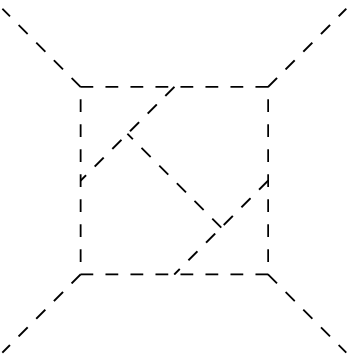}
$}}
\begin{definition}
Let $\calB$ be the space of formal linear
combinations of ``free Jacobi diagrams'' (Jacobi diagrams as before,
but with unmarked univalent ends (``legs'') replacing the circle; see
an example on the right), modulo the $AS$ and $IHX$ relations of
before. Let $\chi:\calB\to\calA$ be the symmetrization map which maps a
$k$-legged free Jacobi diagram to the average of the $k!$ ways of
planting these legs along a circle.
\end{definition}

\parpic(28mm,24mm)[r]{\raisebox{-2mm}{$\xymatrix{
  \calB \ar[r]^\chi \ar[d]^{\calT_\frakg} & \calA \ar[d]^{\calT_\frakg} \\
  \calS(\frakg) \ar[r]^{\chi_\frakg} & \calU(\frakg)
}$}}
\begin{theorem} \label{thm:PBW} (diagrammatic PBW,
\cite{Kontsevich:Vassiliev, Bar-Natan:OnVassiliev}) $\chi$ is an
isomorphism of vector spaces. Furthermore, fixing a metrized $\frakg$ 
there is a commutative square as on the right.
\end{theorem}

Note that $\calB$ can be graded (by half the number of vertices in a Jacobi
diagram) and that $\chi$ respects degrees so it extends to an isomorphism
$\chi:\hat\calB\to\hat\calA$ of graded completions.

\section{Proofs of the Fundamental Theorem} \label{sec:Fundamental}

The heart of all known proofs of Theorem~\ref{thm:Fundamental} is always a
construction of a ``universal finite type invariant'' (see below); it is 
simple to show that the existence of a universal finite type invariant is
equivalent to Theorem~\ref{thm:Fundamental}.

\begin{definition} A universal finite type invariant is a map
$Z:\{\text{knots}\}\to\hat\calA^r$ whose extension to singular knots
satisfies $Z(K)=D+(\text{higher degrees})$ whenever a singular knot $K$
and a chord diagram $D$ are related as in
Section~\ref{subsec:ChordDiagrams}.
\end{definition}

\subsection{The Kontsevich Integral} \label{sec:KontsevichIntegral}
The first construction of a universal
finite type invariant was given by Kontsevich~\cite{Kontsevich:Vassiliev}
(see also~\cite{Bar-Natan:OnVassiliev, ChmutovDuzhin:KontsevichIntegral}).
It is known as ``the Kontsevich Integral'' and up to a normalization factor
it is given by
\[ Z_1(K)=\sum_{m=0}^\infty \frac{1}{(2\pi i)^m}
  \hspace{-4mm}
  \sumint_{\stackrel{t_1<\ldots<t_m}{P=\{(z_i,z'_i)\}}}
  \hspace{-4mm}
  (-1)^{\#P_{\downarrow}}D_P
  \bigwedge_{i=1}^{m}\frac{dz_i-dz'_i}{z_i-z'_i},
\]
where the relationship between the knot $K$, the pairing $P$, the real
variables $t_i$, the complex variables $z_i$ and $z'_i$ and the chord
diagram $D_P$ is summarized by the figure
\[ \eps{width=3in}{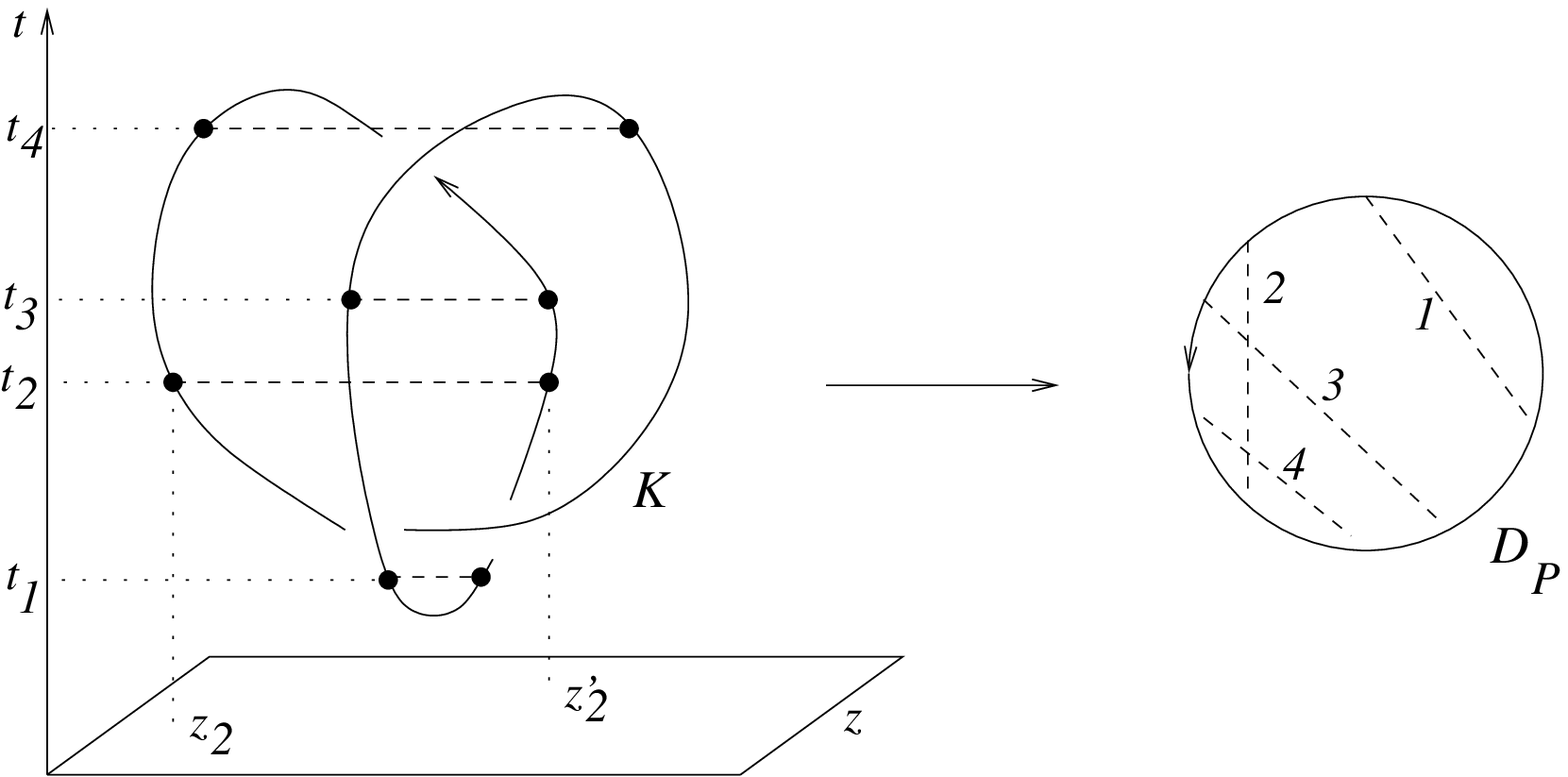}. \]

The Kontsevich Integral arises from studying the holonomy of the
Knizhnik-Zamolodchikov equation of conformal field
theory~\cite{KnizhnikZamolodchikov:CurrentAlgebra}. When evaluating
$Z_1$ one encounters multiple $\zeta$-numbers~\cite{LeMurakami:HOMFLY}
in a substantial way, and the proof that the end result is rational is
quite involved~\cite{LeMurakami:Universal} and relies on deep results
about associators and quasitriangular Quasi-Hopf
algebras~\cite{Drinfeld:QuasiHopf, Drinfeld:GalQQ}.  Employing the same
techniques, in~\cite{LeMurakami:Universal} it is also shown that the
composition of $W_{\frakg,R}\circ Z_1$ precisely reproduces the
Reshetikhin Turaev invariants~\cite{ReshetikhinTuraev:Ribbon}.

\subsection{Perturbative Chern-Simons-Witten theory and configuration
space integrals} \label{sec:CS} Historically the first approach to the
construction of a universal finite type invariant was to use
perturbation theory with the Chern-Simons-Witten topological quantum
field theory; this is also how the relationship with Lie algebras first
arose~\cite{Bar-Natan:Thesis}. But taming the integrals involved turned
out to be difficult and working constructions using this approach
appeared only a bit later~\cite{BottTaubes:SelfLinking,
Thurston:IntegralExpressions, AltschulerFreidel:AllOrders}.

In short, one writes a perturbative expansion for the large $k$
asymptotics of the Chern-Simons-Witten path integral for some metrized
Lie algebra $\frakg$ with a Wilson loop in some representation $R$ of
$\frakg$,
\[
  \int_{\frakg\mbox{\scriptsize -connections}}
  \hspace{-42pt} \calD A\,
  \mbox{\it tr}_R\mbox{\it hol}_K(A)\exp\left[{\scriptstyle
    \frac{ik}{4\pi}\int\limits_{\bbR^3}
    \mbox{tr}\left(A\wedge dA+\frac{2}{3}A\wedge A\wedge A\right)
  }\right].
\]
The result is of the form
\[ 
  \sum_{\parbox{0.7in}{\centering\scriptsize $D$: Feynman diagram}}
    W_{\frakg(D),R}\sumint\calE(D)
\]
where $\calE(D)$ is a very messy integral expression and the diagrams
$D$ as well as the weights $W_{\frakg(D),R}$ are as in
Section~\ref{subsec:Jacobi}. Replacing $W_{\frakg(D),R}$ by simply $D$ in
the above formula we get an expression with values in $\hat\calA$:
\[ Z_2(K):=\sum_D D\sumint\calE(D)\in\hat\calA. \]

For formal reasons $Z_2(K)$ ought to be a universal finite type invariant,
and after much work taming the $\calE(D)$ factors and after multiplying by a
further framing-dependent renormalization term $Z^{\text{anomaly}}$, the
result is indeed a universal finite type invariant.

Upon further inspection the $\calE(D)$ factors can be reinterpreted as
integrals of certain spherical volume forms on certain (compactified)
configuration spaces~\cite{BottTaubes:SelfLinking}. These integrals can
be further interpreted as counting certain ``tinker toy constructions''
built on top of $K$~\cite{Thurston:IntegralExpressions}. The latter
viewpoint makes the construction of $Z_2$ visually
appealing~\cite{Bar-Natan:AstrologyToTopology}, but there is no
satisfactory writeup of this perspective yet.

We note that the precise form of the renormalization term
$Z^{\text{anomaly}}$ remains an open problem. An appealing conjecture is
that $Z^{\text{anomaly}}=\exp\frac12\isolatedchord$. If this is true then
$Z_2=Z_1$ \cite{Poirier:LimitConfigurationSpace}; but the conjecture is
only verified up to degree 6 \cite{Lescop:ForLinks} (there's also an
unconfirmed verification to all orders~\cite{Yang:DegreeTheory}).

The most important open problem about perturbative Chern-Simons-Witten
theory is not directly about finite type invariants, but it is
nevertheless worthwhile to recall it here:

\begin{problem} Does the perturbative expansion of the Chern-Simons-Witten
theory converge (or is asymptotic to) the exact solution due to
Witten~\cite{Witten:Jones} and
Reshetikhin-Turaev~\cite{ReshetikhinTuraev:Ribbon} when the parameter $k$
converges to infinity?
\end{problem}

\subsection{Associators and trivalent graphs} There is also an entirely
algebraic approach for the construction of a universal finite type
invariant $Z_3$. The idea is to find some algebraic context within
which knot theory is finitely presented --- i.e., presented by finitely
many generators subject to finitely many relations. If the algebraic
context at hand is compatible with the definitions of finite type
invariants and of chord diagrams, one may hope to define $Z_3$ by
defining it on the generators in such a way that the relations are
satisfied. Thus the problem of defining $Z_3$ is reduced to finding
finitely many elements of $\calA$-like spaces which solve certain
finitely many equations.

A concrete realization of this idea is in~\cite{LeMurakami:Universal,
Bar-Natan:NAT} (following ideas from~\cite{Drinfeld:QuasiHopf,
Drinfeld:GalQQ} on quasitriangular Quasi-Hopf algebras). The relevant
``algebraic context'' is a category with
certain extra operations, and within it, knot theory is generated by
just two elements, the braiding $\overcrossing$ and the re-association
$\Associator$. Thus to define $Z_3$ it is enough to find
$R=Z_3(\overcrossing)$ and ``an associator'' $\Phi=Z_3(\Associator)$ which
satisfy certain normalization conditions as well as the pentagon and
hexagon equations
\[
  \Phi^{123}\cdot(1\Delta 1)(\Phi)\cdot\Phi^{234}
  =(\Delta 1 1)(\Phi)\cdot
    (1 1\Delta)(\Phi),
\]
\[
  (\Delta 1)(R^{\pm})
  = \Phi^{123} (R^{\pm})^{23}(\Phi^{-1})^{132}
    (R^{\pm})^{13}\Phi^{312}.
\]

As it turns out the solution for $R$ is easy and nearly canonical. But
finding an associator $\Phi$ is a lot harder. There is a closed form
integral expression $\Phi^{KZ}$ due to~\cite{Drinfeld:QuasiHopf} but
one encounters the same not-too-well-understood multiple $\zeta$
numbers of Section~\ref{sec:KontsevichIntegral}. There is a rather
complicated iterative procedure for finding an
associator~\cite{Drinfeld:GalQQ, Bar-Natan:NAT, Bar-Natan:Associators}.
On a computer it had been used to find an associator up to degree $7$.
There is also closed form associator that works only with the Lie
super-algebra $gl(1|1)$~\cite{Lieberum:gl11}. But it remains an open
problem to find a closed form formula for a rational associator (existence
by~\cite{Drinfeld:GalQQ, Bar-Natan:Associators}).

On the positive side we should note that the end result, the invariant
$Z_3$, is independent of the choice of $\Phi$ and that
$Z_3=Z_1$~\cite{LeMurakami:Universal}.

There is an alternative (more symmetric and intrinsicly 3-dimensional,
but less well documented) description of the theory of associators in
terms of knotted trivalent
graphs~\cite{Bar-NatanThurston:AlgebraicStructures, Thurston:Shadow}.
There ought to be a perturbative invariant associated with knotted
trivalent graphs in the spirit of Section~\ref{sec:CS} and such an
invariant should lead to a simple proof that $Z_2=Z_3=Z_1$. But the
$\calE(D)$ factors remain untamed in this case.

\subsection{Step by step integration} \label{subsec:StepByStep}
The last approach for proving the Fundamental Theorem is the most
natural and historically the first. But here it is last because it is
yet to lead to an actual proof. A weight system $W:\calA^r_m\to\bbQ$ is
an invariant of $m$-singular knots. We want to show that it is the
$m$th derivative of an invariant $V$ of non-singular knots. It is
natural to try to integrate $W$ step by step, first finding an
invariant $V^{m-1}$ of $(m-1)$-singular knots whose derivative in the
sense of~\eqref{eq:doublepoint} is $W$, then an invariant $V^{m-2}$ of
$(m-2)$-singular knots whose derivative is $V^{m-1}$, and so on all the
way up to an invariant $V^0=V$ whose $m$th derivative will then be $W$.
If proven, the following conjecture would imply that such an inductive
procedure can be made to work:

\begin{conjecture} \cite{Hutchings:SingularBraids}
If $V^r$ is a once-integrable invariant of $r$-singular
knots then it is also twice integrable. That is, if there is an invariant
$V^{r-1}$ of $(r-1)$-singular knots whose derivative is $V^r$, then there
is an invariant $V^{r-2}$ of $(r-2)$-singular knots whose second
derivative is $V^r$.
\end{conjecture}

In~\cite{Hutchings:SingularBraids} Hutchings reduced this conjecture to a
certain appealing topological statement and further to a certain
combinatorial-algebraic statement about the vanishing of a certain homology
group $H^1$ which is probably related to Kontsevich's graph homology
complex~\cite{Kontsevich:FeynmanDiagrams} (Kontsevich's $H^0$ is $\calA$,
so this is all in the spirit of many deformation theory problems where
$H^0$ enumerates infinitesimal deformations and $H^1$ is the obstruction to
globalization). Hutchings~\cite{Hutchings:SingularBraids} was also able to
prove the vanishing of $H^1$ (and hence reprove the Fundamental Theorem)
in the simpler case of braids. But no further progress has been made along
these lines since~\cite{Hutchings:SingularBraids}.

\section{Some further directions} \label{sec:FurtherDirections}

We would like to touch on a number of significant further directions in the
theory of finite type invariants. We will only say a few words on each of
those and refer the reader to the literature for further information.

\subsection{The original ``Vassiliev'' perspective} V.A.~Vassiliev came
to the study of finite type knot invariants by studying the infinite
dimensional space of all immersions of a circle into $\bbR^3$ and the
topology of the ``discriminant'', the locus of all singular immersions
within the latter space~\cite{Vassiliev:CohKnot, Vassiliev:Book}.
Vassiliev studied the topology of the complement of the discriminant
(the space of embeddings) using a certain spectral sequence and found
that certain terms in it correspond to finite type invariants. This
later got related to the Goodwillie calculus and back to the
configuration spaces of Section~\ref{sec:CS}.
See~\cite{Volic:TaylorTowers}.

\subsection{Interdependent modifications} The standard
definition of finite type invariants is based on modifying a knot by
replacing over (or under) crossings with under (or over) crossings.
In~\cite{Goussarov:Modifications} Goussarov generalized this by allowing
arbitrary modifications done to a knot --- just take any segment of the
knot and move it anywhere else in space. The resulting new ``finite type''
theory turns out to be equivalent to the old one though with a factor of
$2$ applied to the grading (so an ``old'' type $m$ invariant is a ``new''
type $2m$ invariant and vice versa). See also~\cite{Bar-Natan:Bracelets,
Conant:OnGoussarov}.

\subsection{$n$-equivalence, commutators and claspers} While little is
known about the overall power of finite type invariants, much is known
about the power of type $n$ invariants for any given $n$.
Goussarov~\cite{Goussarov:nEquivalence} defined the notion of
$n$-equivalence: two knots are said to be ``$n$-equivalent'' if all
their type $n$ invariants are the same. This equivalence relation is
well understood both in terms of commutator subgroups of the pure braid
group~\cite{Stanford:ModuloPureBraids, NgStanford:GusarovGroup} and in
terms of Habiro's calculus of surgery over
``claspers''~\cite{Habiro:Claspers} (the latter calculus also gives a
topological explanation for the appearance of Jacobi diagrams as in
Section~\ref{subsec:Jacobi}). In particular, already
Goussarov~\cite{Goussarov:nEquivalence} shows that the set of
equivalence classes of knots modulo $n$-equivalence is a finitely
generated Abelian group $G_n$ under the operation of connected sum, and the
rank of that group is equal to the dimension of the space of type $n$
invariants.

Ng~\cite{Ng:Ribbon} has shown that ribbon knots generate an index $2$
subgroup of $G_n$.

\subsection{Polynomiality and Gauss sums}
Goussarov~\cite{Goussarov:PresentedByGauss} (see
also~\cite{GoussarovPolyakViro:VirtualKnots}) found an intriguing way
to compute finite type invariants from a Gauss diagram presentation of
a knot, showing in particular that finite type invariants grow as
polynomials in the number of crossings $n$ and can be computed in
polynomial time in $n$ (though actual computer programs are still
missing!).

Gauss diagrams are obtained from knot diagrams in much of
the same way as Chord diagrams are obtained from singular knots, except all
crossings are counted and not just the double points, and certain
over/under and sign information is associated with each crossing/chord so
that the knot diagram can be recovered from its Gauss diagram. In the
example below, we also dashed a subdiagram of the Gauss diagram equivalent
to the chord diagram $\eps{width=5mm}{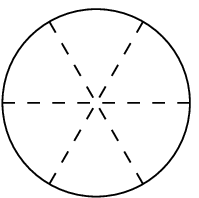}$:
\[ \eps{width=2.4in}{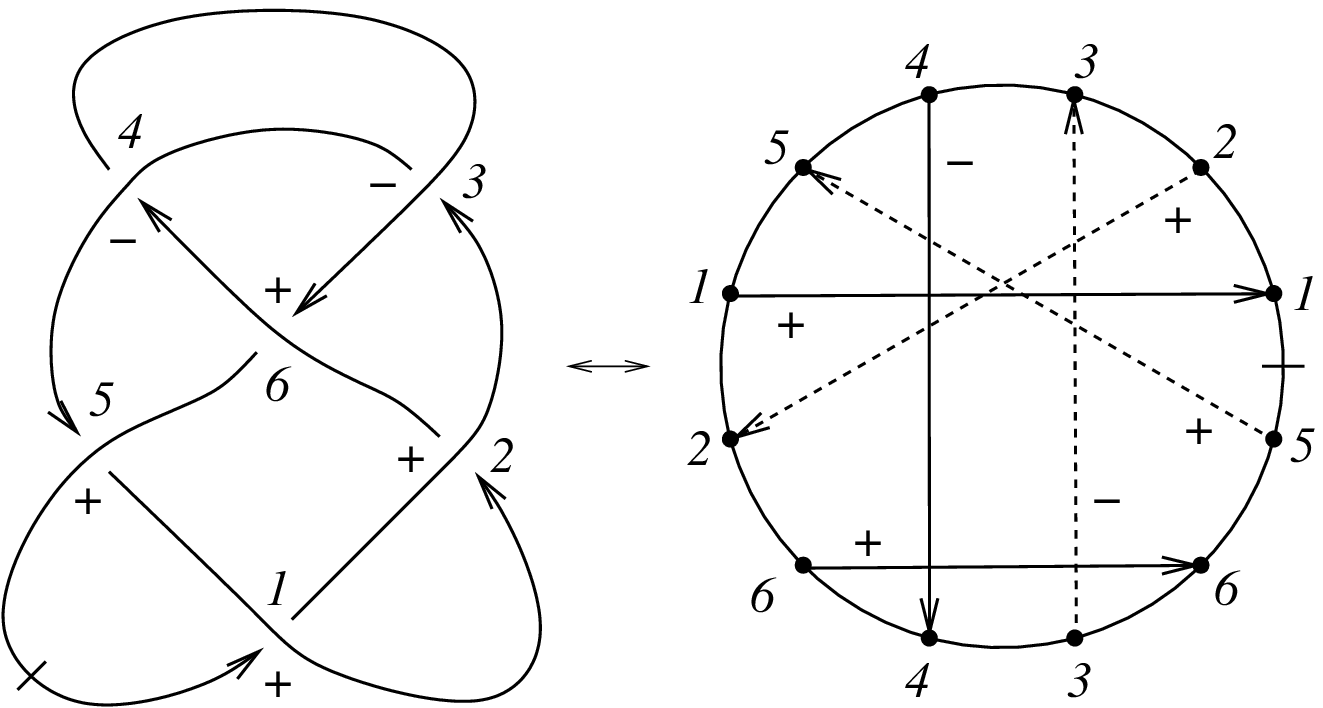} \]

If $G$ is a Gauss diagram and $D$ is a chord diagram we let $\langle
D,G\rangle$ be the number of subdiagrams of $G$ equivalent to $D$, counted
with appropriate signs (to be precise, we also need to base the
diagrams involved and respect the basing).

\begin{theorem} \cite{Goussarov:PresentedByGauss,
GoussarovPolyakViro:VirtualKnots} If $V$ is a type $m$ invariant then there
are finitely many (based) chord diagrams $D_i$ with at most $m$ chords
and rational numbers $\alpha_i$ so that $V(K)=\sum_i\alpha_i\langle
D_i,G\rangle$ whenever $G$ is a Gauss diagram representing a knot $K$.
\end{theorem}

\subsection{Computing the Kontsevich Integral} While the Kontsevich
Integral $Z_1$ is a cornerstone of the theory of finite type invariants,
it has been computed for surprisingly few knots. Even for the unknot the
result is non-trivial:

\begin{theorem} \label{thm:Wheels} (``Wheels'',
\cite{Bar-NatanGaroufalidisRozanskyThurston:WheelsWheeling,
Bar-NatanLeThurston:TwoApplications}) The framed Kontsevich integral of
the unknot, $Z_1^F(\bigcirc)$, expressed in terms of diagrams in
$\hat\calB$, is given by $\Omega=\exp_\udot \sum_{n=1}^\infty
b_{2n}\omega_{2n}$, where the `modified Bernoulli numbers' $b_{2n}$ are
defined by the power series expansion $\sum_{n=0}^\infty b_{2n}x^{2n} =
\frac{1}{2}\log\frac{\sinh x/2}{x/2}$, the `$2n$-wheel' $\omega_{2n}$
is the free Jacobi diagram made of a $2n$-gon with $2n$ legs (so e.g.,
$\omega_6=\eps{width=5mm}{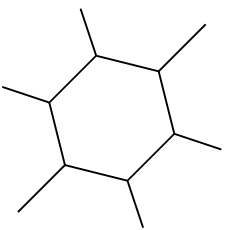}$) and where $\exp_\udot$ means
`exponential in the disjoint union sense'.
\end{theorem}

Closed form formulas have also been given for the Kontsevich integral
of framed unknots, the Hopf link and Hopf
chains~\cite{Bar-NatanLawrence:RationalSurgery} and for torus
knots~\cite{Marche:Computation}.

Theorem~\ref{thm:Wheels} has a companion that utilizes the same element
$\Omega$, the ``wheeling''
theorem~\cite{Bar-NatanGaroufalidisRozanskyThurston:WheelsWheeling,
Bar-NatanLeThurston:TwoApplications}. The wheeling theorem ``upgrades''
the vector space isomorphism $\chi:\calB\to\calA$ to an algebra
isomorphism and is related to the Duflo isomorphism of the theory of
Lie algebras. It is amusing to note that the wheeling theorem (and hence
Duflo's theorem in the metrized case) follows using finite type techniques
from the ``$1+1=2$ on an abacus'' identity
\[ \eps{width=2.6in}{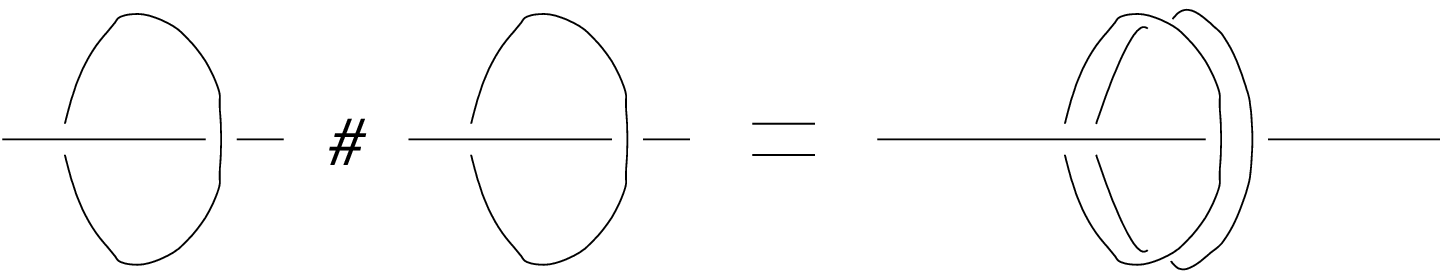}. \]

\subsection{Taming the Kontsevich Integral} While explicit calculations
are rare, there is a nice structure theorem for the values of the
Kontsevich integral, saying that for a knot $K$ and up to any fixed
number of loops in the Jacobi diagrams, $\chi^{-1}Z_1(K)$ can be
described by finitely many rational functions (with denominators powers
of the Alexander polynomial) which dictate the placement of the legs. This
structure theorem was conjectured in~\cite{Rozansky:RationalityConjecture},
proven in~\cite{Kricker:RationalityConjecture} and partially generalized to
links in~\cite{GaroufalidisKricker:NonCommutativeInvariant}.

\subsection{The Rozansky-Witten theory} One way to construct linear
functionals on $\calA$ (and hence finite type invariants) is using Lie
algebras and representations as in Section~\ref{subsec:Jacobi}; much of
our insight about $\calA$ comes this way. But there is another
construction for such functional (and hence invariants), due to
Rozansky and Witten~\cite{RozanskyWitten:HyperKahler}, using
contractions of curvature tensors on hyper-K\"ahler manifolds. Very
little is known about the Rozansky-Witten approach; in particular, it
is not known if it is stronger or weaker than the Lie algebraic
approach. For an application of the Rozansky-Witten theory back to
hyper-K\"ahler geometry check~\cite{HitchinSawon:HyperKahler} and for a
unification of the Rozansky-Witten approach with the Lie algebraic
approach (albeit at a categorical level)
check~\cite{RobertsWillerton:InPreparation}.

\subsection{The Melvin-Morton conjecture and the volume conjecture} The
Melvin-Morton conjecture (stated~\cite{MelvinMorton:Coloured},
proven~\cite{Bar-NatanGaroufalidis:MMR}) says that the Alexander
polynomial can be read off certain coefficients of the coloured Jones
polynomial. The Kashaev-Murakami-Murakami volume conjecture
(stated~\cite{Kashaev:HyperbolicVolume,
MurakamiMurakami:SimplicialVolume}, unproven) says that a certain
asymptotic growth rate of the coloured Jones polynomial is the
hyperbolic volume of the knot complement. 

Both conjectures are not directly about finite type invariants but both
have ramifications to the theory of finite type invariants. The
Melvin-Morton conjecture was first proven using finite type invariants
and several later proofs and generalizations
(see~\cite{Bar-Natan:VasBib}) also involve finite type invariants. The
volume conjecture would imply, in particular, that the hyperbolic
volume of a knot complement can be read from that knot's finite type
invariants, and hence finite type invariants would be at least as
strong as the volume invariant.

A particularly noteworthy result and direction for further research is
Gukov's~\cite{Gukov:APolynomial} recent unification of these two
conjectures under the Chern-Simons umbrella (along with some relations
to three dimensional quantum gravity).

\section{Beyond knots} \label{sec:BeyondKnots} For the lack of space we
have restricted ourselves here to a discussion of finite type
invariants of knots. But the basic ``differentiation'' idea of
Section~\ref{sec:intro} calls for generalization, and indeed it has been
generalized extensively. We will only make a few quick comments.

Finite type invariants of homotopy links (links where each component is
allowed to move across itself freely) and of braids are extremely well
behaved. They separate, they all come from Lie algebraic
constructions~\cite{Bar-Natan:Homotopy, Lin:Expansions,
Bar-Natan:Braids} and in the case of braids, step by step integration
as in Section~\ref{subsec:StepByStep}
works~\cite{Hutchings:SingularBraids} (for homotopy links the issue was
not studied).

Finite type invariants of 3-manifolds and especially of integral and
rational homology spheres have been studied extensively and the picture
is nearly a complete parallel of the picture for knots. There are
several competing definitions of finite type invariants (due to
\cite{Ohtsuki:IntegralHomology} and then
\cite{Goussarov:PresentedByGauss, Garoufalidis:3ManifoldsI,
GaroufalidisGoussarovPolyak:Clovers} and more), and they all agree up to
regrading. There are weight systems and they are linear functionals on a
space $\calA(\emptyset)$ which is a close cousin of $\calA$ and $\calB$ and
is related to Lie algebras and hyper-K\"ahler manifolds in a similar way.
There is a notion of a ``universal'' invariant, and there are several
constructions (due to \cite{LeMurakamiOhtsuki:Universal, Le:UniversalIHS}
and then \cite{Bar-NatanGaroufalidisRozanskyThurston:Aarhus,
KuperbergThurston:CutAndPaste}), they all agree or are conjectured to
agree, and they are related to the Chern-Simons-Witten theory.

Finite type invariants were studied for several other types of
topological objects, including knots within other manifolds, higher
dimensional knots, virtual knots, plane curves and doodles and more.
See~\cite{Bar-Natan:VasBib}.

\end{multicols}

\end{document}